\documentclass[aps,prl,reprint,superscriptaddress,amsmath,amssymb,showpacs]{revtex4-1}
\usepackage{graphicx}
\usepackage{dcolumn}
\usepackage{bm}
\usepackage{hyperref}
\usepackage[mathlines]{lineno}
\usepackage{epstopdf} 
\usepackage{xcolor}

\DeclareMathOperator{\erf}{\operatorname{erf}}
\newcommand{\eps}{\varepsilon}
\newcommand{\Iapp}{I_{\text{app}}}
\newcommand\TV[1]{{\color{black}{#1}}}

\begin{document}

\preprint{APS/123-QED}
\title{\TV{Delayed bifurcation phenomena in reaction--diffusion equations:\\persistence of canards and slow passage through Hopf bifurcations}}
\author{Tasso J. Kaper}
\affiliation{Department of Mathematics and Statistics, Boston University, Boston, MA 02215, USA}
\author{Theodore Vo}
\affiliation{Department of Mathematics, Florida State University, Tallahassee, FL 32306, USA}
\date{\today}

\begin{abstract}
\TV{
In the context of a spatially extended model for the electrical activity in a pituitary lactotroph cell line, we establish that two delayed bifurcation phenomena from ODEs ---folded node canards and slow passage through Hopf bifurcations--- persist in the presence of diffusion.  
For canards, the single cell (ODE) model exhibits canard-induced bursting.
Numerical simulations of the PDE reveal rich spatio-temporal canard dynamics, and the transitions between different bursts are mediated by spatio-temporal maximal canards. 
The ODE model also exhibits delayed loss of stability due to slow passage through Hopf bifurcations.
Numerical simulations of the PDE reveal that this delayed stability loss persists in the presence of diffusion. 
To quantify and predict the delayed loss of stability, we show that the Complex Ginzburg-Landau equation exhibits the same property, and derive a formula for the space-time boundary that acts as a buffer curve beyond which the delayed onset of oscillations must occur. 
}
\end{abstract}

\maketitle


\paragraph{Introduction}
\TV{
Canards and maximal canards \cite{Benoit1981,Diener1984,Dumortier1996,Eckhaus1983,Krupa2001} arise ubiquitously in multi-scale ODEs \cite{Bold2003,Brons2006,Carter2017,Desroches2012,Drover2004,Hasan2017,Moehlis2002,Rubin2007}.
They are generic in systems with at least two slow variables \cite{Szmolyan2001,Wexy2005,Wexy2012} and a form of bifurcation delay. 
Maximal canards are the phase and parameter space boundaries between different rhythms, such as spiking and bursting \cite{Brons1991,Burke2012,Burke2016,Desroches2016,Mitry2013,Rotstein2003,Vo2017}.

Delayed Hopf bifurcation (DHB) \cite{Shishkova1973, Neishtadt1987,*Neishtadt1988, Baer1989, Neishtadt1995,Hayes2016} is the other ubiquitous delayed bifurcation phenomenon in ODEs \cite{Ashwin2012, Barreto2011,Erneux1990,Erneux1991,Holden1993,Koper1998, Kuehn2015, Rubin2002, Shorten2000,Strizhak1996}. 
Solutions starting near stable quasi-stationary states (QSS) stay near them for long times after they have become unstable in Hopf bifurcations. 
Substantial delays in the stability loss --and hence in the onset of oscillations-- are observed.
In neuronal models, DHB arises because the intracellular calcium concentration varies slowly and moves the system through Hopf bifurcations \cite{Bertram1995,DelNegro1998,Izhikevich2001,Rinzel1988}.

In this Letter, we first report on the persistence of canards and maximal canards in the presence of diffusion. 
This new PDE phenomenon is reported for a model of electrical activity in a pituitary lactotroph cell line. 
We find that rich spatio-temporal canard dynamics are responsible for complex bursting rhythms on large open regions of parameter space.
In the same PDE, we report on the discovery of delayed loss of stability due to slow passage through Hopf bifurcations.
The delayed stability loss determines important burst diagnostics, including frequency and amplitude.

To supplement the numerics, we show that the Complex Ginzburg-Landau (CGL) equation \cite{Aranson2002,Eckhaus1965,Manneville1990,vanSaarloos1993,Newell1987,Nicolis1995,Walgraef1997} with slowly varying parameter exhibits DHBs and derive a formula to predict when solutions must diverge from the unstable QSS. 
This formula defines a buffer curve in the space-time plane along which the delayed onset of oscillations must occur. 
The analytical buffer curve agrees with the numerics over a wide range of diffusivities.
}

Finally, we simulated spatio-temporal canards in the forced van der Pol equation with diffusion, and spatially inhomogeneous DHBs in the Hodgkin-Huxley \cite{Borgers2016}, Brusselator \cite{Walgraef1997}, FitzHugh-Nagumo \cite{Su1993,*Su1994}, Morris-Lecar \cite{Meier2015}, and Hindmarsh-Rose \cite{Raghavachari1999} PDEs (not shown), showing both new PDE phenomena occur widely.

\paragraph{\TV{Canards in pituitary cell line model}}

\TV{We report on the existence of spatio-temporal canards in a PDE model of electrical activity in a pituitary lactotroph cell line,
\begin{equation}	\label{eq:model}
\begin{split}
C_m V_t &= - \sum I_{\text{ionic}} + I + \Iapp(x) + D V_{xx}, \\ 
\tau_n n_t &= n_{\infty}(V) - n, \\ 
\tau_e e_t &= e_{\infty}(V) - e.
\end{split}	
\end{equation}
The ionic current, $\sum I_{\text{ionic}}$, consists of calcium ($I_{Ca}$), delayed rectifier K$^+$ ($I_K$), A-type K$^+$ ($I_A$), and leak ($I_L$) currents. 
The state variables $(V,n,e)$ are the membrane potential and gating variables for the activation of $I_K$ and inactivation of $I_A$ \cite{Toporikova2008,Vo2012}; $I$ is a baseline current.
The cells, whose positions are denoted by $x$, are coupled via gap junctions with diffusivity $D$. 

The applied current $\Iapp (x)$ is localized. 
Experiments and models have identified the importance of treating portions of the brain as spatially inhomogeneous media, with localized synaptic currents.
For example, EEG data from the auditory cortex in certain primates exhibits spatially localized currents (Fig. 2 \cite{Lakatosetal2005}); the locally generated intracortical synaptic currents have tapered peaks at supragranular and granular sites.

\begin{figure}[h]
\centering
\includegraphics[width=0.975\columnwidth]{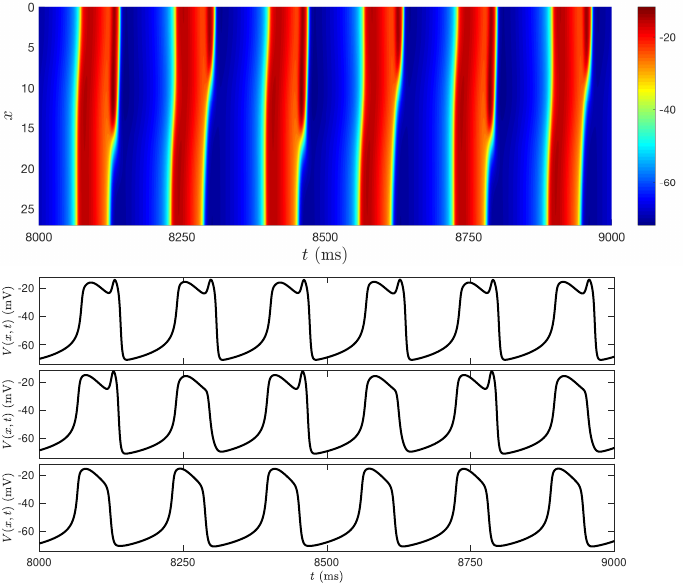}
\put(-244,200){(a)}
\put(-244,109){(b)}
\caption{Simulation of \eqref{eq:model} with $L=50$, $D=1$, $g_K=6.15$~nS, $g_A = 5$~nS, $a = 1$, and $\sigma = 50$. 
(a) Colour map of the voltage, $V(x,t)$, blue (hyperpolarization) and red (depolarization). 
(b) Time series at $x = 1$ ($1^1$ burst; top), $x=11.6$ ($1^1 1^0$ alternator; middle), and $x = 24$ ($1^0$ spike; bottom).
}
\label{fig:alternatortime}
\end{figure}

We performed simulations of \eqref{eq:model} with various $\Iapp(x)$ and zero-flux boundary conditions on $[-L,L]$
\footnote{Numerical simulations were performed using balanced symmetric Strang operator splitting \cite{Speth2013}, with centered finite differences for the Laplacian and fourth-order Runge-Kutta for the time stepping. 
The discretization was chosen fine enough to resolve all modes up to those with $k=\mathcal{O}(\eps^{-3/2})$.
We verified our numerical results independently using Crank-Nicolson.}; results are shown for $x \geq 0$.
A representative simulation with $\Iapp(x) = a \exp ( -\frac{x^2}{4\sigma} )$ is shown in Fig.~\ref{fig:alternatortime}. 
There are three regions. 
In the central region (about $x=0$), the attractor consists of an $L^s = 1^1$ burst (Fig.~\ref{fig:alternatortime}(b), top), with $s=1$ small-amplitude oscillations (SAOs) in the depolarized phase for every $L$ large-amplitude relaxation-type oscillations. 
In the alternator region ($9 \lesssim |x| \lesssim 20$), the attractor alternates periodically between $1^1$ and $1^0$ bursts (Fig.~\ref{fig:alternatortime}(b), middle), where a $1^0$ burst is a spike.   
In the spiking region ($|x| \gtrsim 20$), the attractor consists of $1^0$ spikes (Fig.~\ref{fig:alternatortime}(b), bottom).

The widths of the regions are determined by $\Iapp(x)$. 
For each $x$ in the center, $\Iapp(x)$ is large enough that the $x$-dependent ODE (i.e., \eqref{eq:model} with $D=0$) exhibits $1^1$ canard-induced bursting \cite{Brons2006} attractors over a large domain in $(g_K,g_A)$ parameter space \cite{Vo2012}.
(Increases in $\Iapp$ are equivalent to decreases in $g_K$. The SAOs are folded node canards.)
Here, the state of the PDE lies close to this family of $1^1$ bursts, and converges to it as $D \to 0$.
For each $x$ in the spiking region, $\Iapp(x)$ is small enough that the $x$-dependent ODE exhibits a $1^0$ spiking attractor over a large domain in $(g_K,g_A)$ space.
Here, the state of the PDE lies close to this family of $1^0$ spikes. 
For each $x$ in the alternator region, $\Iapp(x)$ is such that the $x$-dependent ODE exhibits bursts with alternating signatures. 
Here, the state of the PDE is a $1^1 1^0$ rhythm. 

To better understand the PDE dynamics, we compare to the ODE dynamics \cite{Vo2012}. 
With $D=0$, \eqref{eq:model} is an $x$-dependent ODE for each $x$ in $[-L,L]$, and the attracting and repelling invariant slow manifolds, $\mathcal{S}_a(x)$ and $\mathcal{S}_r(x)$, organize the dynamics.
For each $x$, the manifolds were computed using pseudo-arclength continuation \cite{Desroches2008,Desroches2010}.
The intersections of the slow manifolds are {\em maximal canards}, which partition $\mathcal{S}_a(x)$ and $\mathcal{S}_r(x)$ into rotational sectors \cite{Szmolyan2001,Wexy2005}. 
The maximal strong canard, $\gamma_0(x)$, in Fig.~\ref{fig:canardinnerouter}(a) divides the solutions which exhibit local oscillations and those which do not.  

In this simulation, for all $x$ in the bursting region, the steady state lies in the sector with one local oscillation.
It stays close to the family of $\mathcal{S}_a(x)$ until the folded node, then follows the family of $\mathcal{S}_r(x)$ for a long time, after which it transitions to the hyperpolarized state. 
Thus, the $1^1$ canards of the ODEs persist in the bursting region for $D>0$. 
 
For all $x$ in the spiking region, there are no intersections of $\mathcal{S}_a(x)$ and $\mathcal{S}_r(x)$. 
Here, the steady state of the PDE follows the family of $\mathcal{S}_a(x)$ until the folded node, after which it transitions directly to the hyperpolarized state without any SAOs (Fig.~\ref{fig:canardinnerouter}(b)).
Hence, the $1^0$ spikes also persist in the PDE, and can co-exist in the same steady state with the $1^1$ bursts.

\begin{figure}[h]
\centering
\includegraphics[width=\columnwidth]{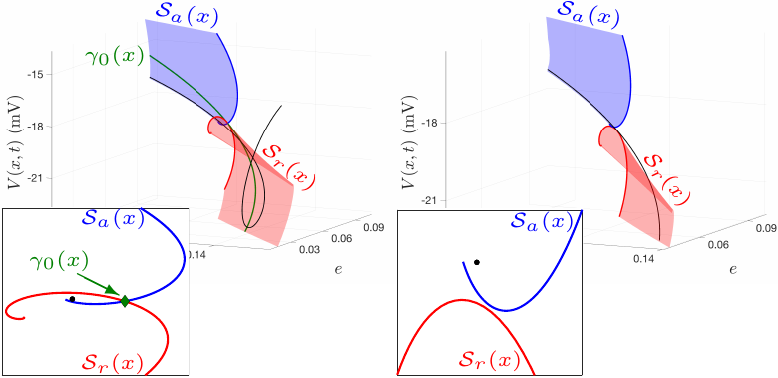}
\put(-250,106){(a)}
\put(-125,106){(b)}
\caption{Steady state of the PDE (black) from Fig.~\ref{fig:alternatortime} compared to $\mathcal{S}_a(x)$ and $\mathcal{S}_r(x)$ at (a) $x=1$ ($1^1$ bursting region), and (b) $x=24$ (spiking region). Insets: cross-section taken through the folded node.}
\label{fig:canardinnerouter}
\end{figure}

We also find that the maximal canards of the ODEs persist in the presence of diffusion. 
The steady state solutions of \eqref{eq:model} exhibit maximal canards in the transition intervals between different regions. 
For example, the steady state solution from Fig.~\ref{fig:alternatortime} exhibits period-2 bursts in the interval between the central and alternator regions. 
For $x$ close to the center, the SAOs of the odd bursts have larger amplitude than those of the even ones (Fig.~\ref{fig:transition}(a)). 
Further from the center, the amplitude of the odd-burst-SAOs is smaller, whilst that of the even ones is larger (Fig.~\ref{fig:transition}(b)--(c)). 
Sufficiently far from the central region, the odd-burst-SAOs are absent and the even-burst-SAOs have maximum amplitude, corresponding to a maximal canard (Fig.~\ref{fig:transition}(d)). 
In this manner, the system transitions in space from the $1^1$ bursting state to the $1^1 1^0$ alternator state. 
Similar attractors and maximal canards are observed for other $D$.

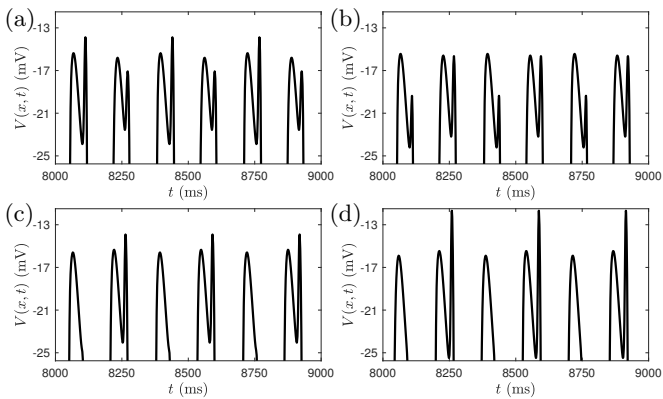
\begin{figure}[h]
\centering
\includegraphics[width=\columnwidth]{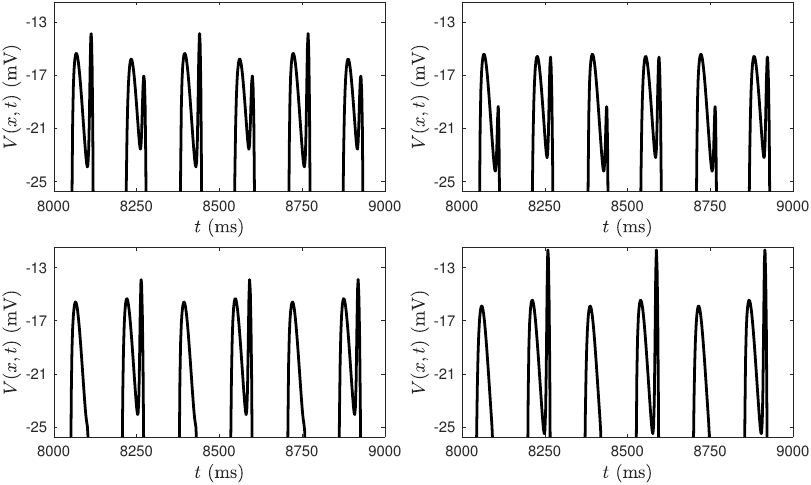}
\put(-249,141){(a)}
\put(-126,141){(b)}
\put(-249,67){(c)}
\put(-126,67){(d)}
\caption{Maximal canards in the transition interval between the central and alternator regions of Fig.~\ref{fig:alternatortime}.  
Slices taken at (a) $x \approx 7.52$, (b) $x \approx 8.94$, (c) $x \approx 10.25$, and (d) $x \approx 13.48$.}
\label{fig:transition}
\end{figure}

In other regions of parameter space, \eqref{eq:model} exhibits $1^s$ canard-induced bursts with $s\geq 1$. 
Fig.~\ref{fig:mmotransition} shows a simulation in which $1^3$ bursts from the central region invade $1^2$ bursts in the outer region via a moving front. 
This transition is also mediated by maximal canards.
The front speed increases with $D$.

\begin{figure}[h]
\centering
\includegraphics[width=\columnwidth]{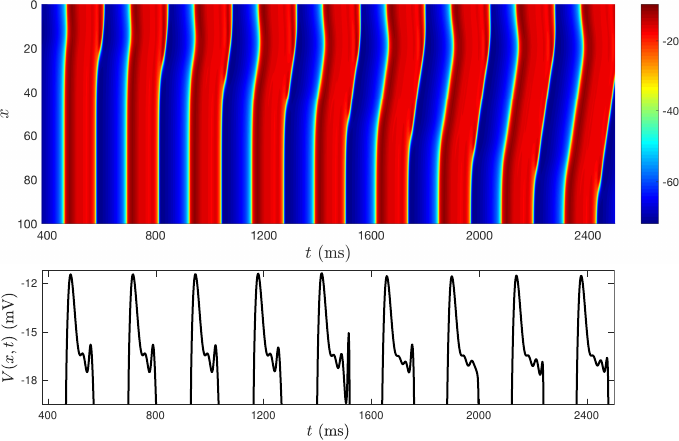}
\put(-249,154){(a)}
\put(-249,60){(b)}
\caption{Solution of \eqref{eq:model} with $L=100$, $D=1$, $g_K=4.35$~nS, $g_A=5$~nS, $a=1$ and $\sigma=50$. 
The $1^3$ bursts from the central region invade the $1^2$ region.
(a) Colour map of $V(x,t)$ and (b) time series at $x = 48$.
}
\label{fig:mmotransition}
\end{figure}

}

\paragraph{DHB in \eqref{eq:model}}

\TV{
Next, we show that the PDE \eqref{eq:model} exhibits delayed loss of stability due to slow passage through Hopf bifurcations (DHB) when the baseline current is slowly varying, $I(t) = I_0 -\eps t$. 
As $\eps \to 0$, \eqref{eq:model} with $D=0$ possesses depolarized QSS (Fig.~\ref{fig:pitslowpass}; red) given by 
\begin{equation}	\label{eq:HHQSS}
\begin{split}
I &= I_{Ca}+I_K+I_A + I_L - \Iapp(x),
\end{split}
\end{equation}
where $n = n_{\infty}(V)$ and $e = e_{\infty}(V)$.
The stable and unstable parts of the QSS are separated by a curve, $\mathcal{H}$, of subcritical HB (Fig.~\ref{fig:pitslowpass}; black).
For $I$ to the right (left) of $\mathcal{H}$, the QSS are stable (unstable, resp.). 

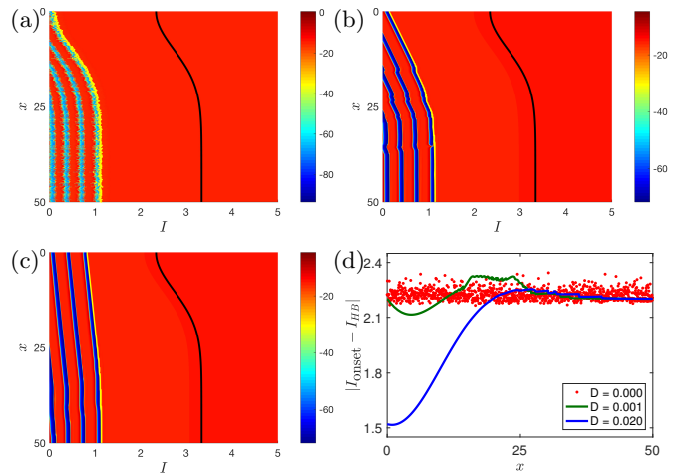
\begin{figure}[h]
\centering
\includegraphics[width=\columnwidth]{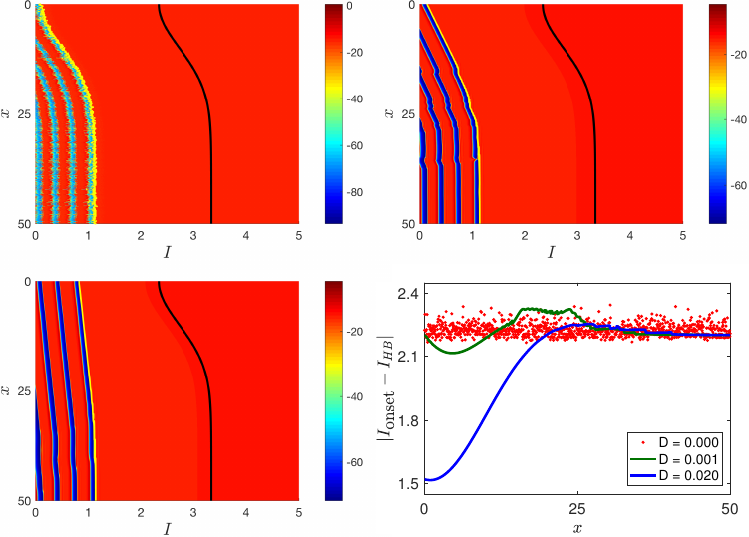}
\put(-249,169){(a)}
\put(-127,169){(b)}
\put(-249,78){(c)}
\put(-127,78){(d)}
\caption{DHB in \eqref{eq:model} with $g_K=4$~nS, $g_A=5$~nS, $a=1$, and $\sigma=50$. (a)--(c) Voltage and (d) measurement of bifurcation delay. Solutions stay close to the QSS (red background) beyond (i.e., left of) $\mathcal{H}$ and transition to bursts at the (yellow) onset curve. (a) $D=0$, (b) $D = 0.001$, and (c) $D = 0.02$.}
\label{fig:pitslowpass}
\end{figure}

With $0< \eps \ll 1$, solutions pass through $\mathcal{H}$ due to the slow decrease in $I$, and there is delayed loss of stability. 
The system stays close to its QSS past $\mathcal{H}$ with canard-induced bursting oscillations (spatially homogeneous or inhomogeneous) setting in at a significantly later time.

\paragraph{Spatially homogeneous DHB}
For homogeneous $\Iapp$, solutions of \eqref{eq:model} with initial conditions sufficiently far from $\mathcal{H}$ (i.e., $|I_0-I_{\text{HB}}|$ large enough) escape the neighbourhood of the QSS \footnote{We classified solutions as having escaped the QSS when their distance from the QSS exceeded $\sqrt{\eps}$.} and oscillate uniformly after a substantial delay. 
The uniformity (for any $D$) is due to the homogeneity of both $\Iapp$ and the frequency along $\mathcal{H}$.
Onset occurs at the value of $I$ predicted from the ODE (not shown).

\paragraph{Spatially inhomogeneous DHB}
Inhomogeneous $\Iapp(x)$ induce inhomogeneous delayed stability loss in \eqref{eq:model}. 
The delay in the stability loss tends to be shortest where $|\Iapp(x)|$ takes its maximum value, and tends to lengthen as the value of $|\Iapp(x)|$ decreases (Fig.~\ref{fig:pitslowpass}). 

For a Gaussian source, we measured the distance, $\left| I_{\rm onset} - I_{\rm HB} \right|$, that solutions stayed close to the QSS past $\mathcal{H}$ (Fig. \ref{fig:pitslowpass}(d)). 
For $D=0$, the delay is almost spatially uniform (red markers). The minor variations are due to the numerical sensitivity associated with using initial value solvers to follow unstable QSS -- more precise measurements can be made by using boundary value solvers. 
For small $D$ (Fig.~\ref{fig:pitslowpass}(b)), the instability that first sets in at $x=0$ spreads locally, as reflected in the minimum in the (green) delay curve.
There is no long-range effect since the green curve lies close to the red for $|x| \gtrsim 25$. 
For larger $D$ (Fig.~\ref{fig:pitslowpass}(c)), there is also no long-range effect. 
However, near the center, the delay duration is much shorter than in the ODE (blue curve). 
Therefore, the PDE \eqref{eq:model} exhibits rich DHB dynamics.
}

\paragraph{DHB in CGL}
To quantify and predict DHB in PDEs, we analyze slow passage through HB in the CGL equation with source term, $\Iapp(x)$, and slowly increasing (real) parameter, $\mu$,
\begin{equation}	\label{eq:CGL}
\begin{split}
A_t &= \left( \mu+i \omega_0 \right) A + \eps D A_{xx} + \sqrt{\eps} \Iapp(x)- \alpha \left| A \right|^2 A, \\
\mu_t &= \eps.
\end{split}
\end{equation}
Here, $A$ is complex, $\omega_0$ is the linear frequency, $\alpha = 1 + i \alpha_i$ is related to the nonlinear frequency, $D = \beta_r + i\beta_i$ is related to the linear dispersion coefficient \cite{Kuramoto1984,Stich2004}, and $0<\eps \ll 1$ measures the timescale separation. 
\TV{We present results for $\omega_0>0$; we find similar results for $\omega_0<0$.}
Simulations of \eqref{eq:CGL} were performed with zero-flux boundary conditions on $[-L,L]$ \footnote{On sufficiently large domains, we observe similar results for the onset of oscillations with Dirichlet conditions.}.
Unless stated otherwise, $\eps = 0.01$, $\omega_0 = 0.5$, $\alpha_i = 0.6$, $\beta_r = 1$, and $\beta_i = 0$. 

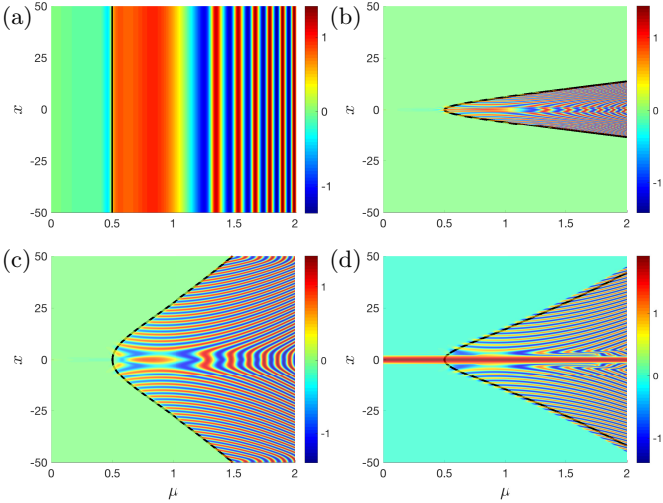
\begin{figure}[htb]  
\centering
\includegraphics[width=\columnwidth]{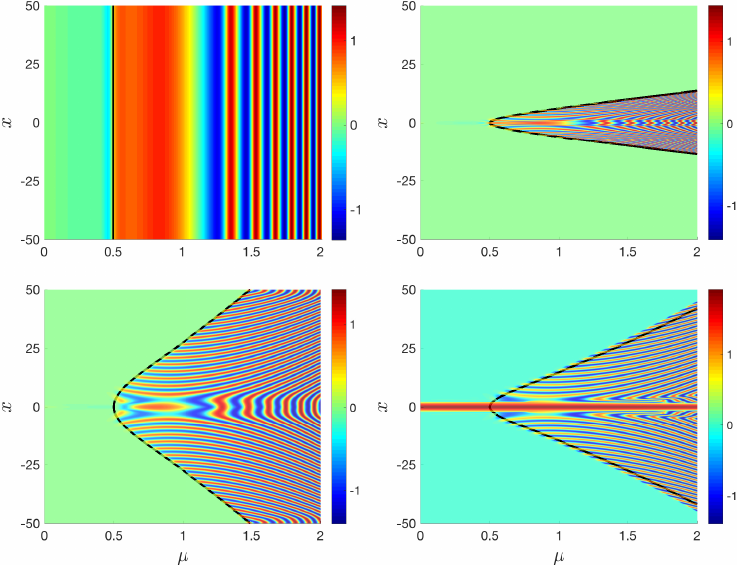}
\put(-250,180){(a)}
\put(-126,180){(b)}
\put(-250,88){(c)}
\put(-126,88){(d)}
\caption{DHB in \eqref{eq:CGL}. (a) Spatially homogeneous source ($a=1$). 
(b)--(d) Gaussian sources ($\sigma = 0.25$) with $(\beta_r,\beta_i,a) = (0,0,1), (3,1,1)$, and $(1,0,100)$.
$\mu(0) = \mu_0=-1$, left of $\mathcal{H}$. 
The buffer curves \eqref{eq:bufferman} are the black curves.
Away from the center, diffusion induces earlier onset of oscillations compared to the diffusionless case, widening the onset curve.  
(Green background indicates QSS.)
}
\label{fig:CGLprofiles}
\end{figure}

In the limit $\eps \to 0$, \eqref{eq:CGL} has a supercritical HB, at $\mu=0$ with frequency $\omega_0$.
The state $A \equiv 0$ is stable (unstable) if $\mu<0$ ($\mu>0$). 
With $0<\eps \ll 1$, the slow increase of $\mu$ causes \eqref{eq:CGL} to pass through the curve, $\mathcal{H}$, of HB (Fig. \ref{fig:CGLprofiles}).

For all sources and with $\mu_0 <-\omega_0$, the system exhibits DHB, i.e., the solution remains close to the QSS,
\begin{equation} \label{eq:AQSS}
\begin{split} 
A = -\tfrac{\Iapp(x)}{\mu + i \omega_0} &\sqrt{\eps} + \left( \tfrac{\Iapp(x) + D (\mu + i \omega_0) \Iapp''(x)}{(\mu + i \omega_0)^3} \right. \\
& \left. - \tfrac{(1+i\alpha) \Iapp^3(x)}{(\mu + i \omega_0)^2(\mu^2 + \omega_0^2)} \right) \eps^{3/2} + \mathcal{O}(\eps^{5/2}),
\end{split}
\end{equation}
well past $\mathcal{H}$ (into the regime where the QSS is repelling). 

For homogeneous sources ($\Iapp(x) = a$), which arise in the problem of spatially-resonant traveling wave forcing of the CGL equation for $1:m$ resonances \cite{Rudiger2007,Ma2011}, the transition from QSS to temporal oscillations occurs uniformly with a significant delay beyond $\mathcal{H}$.
The delay is the maximal delay time predicted by the theory of DHB in ODEs \cite{Shishkova1973,Neishtadt1987,Neishtadt1988,Neishtadt1995,Baer1989,Su1993,Hayes2016}, namely at $\mu = 0.5 = \omega_0$ (Fig.~\ref{fig:CGLprofiles}(a)).

For inhomogeneous sources, the instability first occurs at time $\mu = \omega_0$, where $\left| \Iapp(x) \right|$ is maximal.
For Gaussians, the oscillations set in at $(\mu,x)=(0.5,0)$. 
For larger $|x|$, the instability occurs for later $\mu$, making the delay spatially inhomogeneous (Fig. \ref{fig:CGLprofiles}(b)--(d)). 
This defines the \emph{onset curve} in the $(\mu,x)$ plane.

\paragraph{Buffer curves}
There is a {\em buffer curve} in the $(\mu,x)$ plane past which solutions of \eqref{eq:CGL} with any $\mu_0 < -\omega_0$ cannot stay near the repelling QSS. 
For general $\Iapp(x)$, the buffer curves are derived by Fourier analysis of the linearization about the QSS,
$\eps \hat{A}_{\mu} = \left( \mu +i \omega_0 - \eps D k^2 \right) \hat{A} + \sqrt{\eps} \hat{I}_{\text{app}}(k)$,
where $\hat{A}$ is the Fourier transform of $A$. 
With $\hat{A}(\mu_0,k) = \hat{A}_0$, the solution has homogeneous part,
\begin{equation*} 
\begin{split}
\hat{A}_{\text{hom}} = \hat{A}_0 \exp \left( \tfrac{1}{2\eps} (\mu-\mu_0)(\mu+\mu_0+2i \omega_0-2\eps D k^2) \right),
\end{split}
\end{equation*}
and inhomogeneous part
\begin{equation*} 
\begin{split}
\hat{A}_{\text{inhom}} = \sqrt{\tfrac{\pi}{2}}& \hat{I}_{\text{app}}(k) \exp \left( \tfrac{(\mu+i \omega_0-\eps D k^2)^2}{2\eps} \right) \times \\
& \left[ \erf \left( \tfrac{\mu+i \omega_0 - \eps D k^2}{\sqrt{2\eps}} \right) - \erf \left( \tfrac{\mu_0+i \omega_0 - \eps D k^2}{\sqrt{2\eps}} \right) \right],
\end{split}
\end{equation*}
\footnote{Here, $\erf$ is the error function, 
$ \erf(z) = \frac{2}{\sqrt{\pi}} \int_0^z e^{-t^2}\, dt$.}.
Fourier inversion shows that $A_{\text{hom}}$ is
\begin{equation} \label{eq:hom} 
\tfrac{A(x,\mu_0)}{\sqrt{2D(\mu-\mu_0)}} \ast \exp \left( \tfrac{(\mu-\mu_0)(\mu+\mu_0+2i \omega_0)}{2\eps} - \tfrac{x^2}{4D(\mu-\mu_0)} \right),
\end{equation}
where the asterisk denotes the convolution. Similarly, $A_{\text{inhom}}$ is asymptotically proportional to 
\begin{equation} \label{eq:inhom} 
\Iapp(x) \ast \exp \left( \tfrac{(\mu+i \omega_0)^2}{2\eps} - \tfrac{x^2}{4 D (\mu+i \omega_0)} \right),
\end{equation}
neglecting $\mathcal{O}\left( \eps^2 D^2 k^4 \right)$ and higher order terms. 

The buffer curves are determined by \eqref{eq:inhom}. 
For initial conditions $\mu_0 <- \omega_0$, both \eqref{eq:hom} and \eqref{eq:inhom} remain exponentially small (in $\eps$) until at least $\mu = \omega_0$. Moreover, \eqref{eq:inhom} grows exponentially before \eqref{eq:hom} does. 
Thus, the curve in the $(\mu,x)$ plane along which \eqref{eq:inhom} first begins to grow exponentially is the buffer curve past which solutions cannot remain near the repelling QSS.
This buffer curve is the PDE analog of the buffer point for DHB in ODEs \cite{Neishtadt1987,Neishtadt1988}.

\paragraph{Buffer curves for Gaussian sources}
Evaluation of \eqref{eq:inhom} for a Gaussian source shows that the transition from decay to growth occurs along the following buffer curve:
\begin{equation} \label{eq:bufferman}
\mu^2 = \omega_0^2 + \frac{\eps x^2 \left( \sigma+ \mu \beta_r - \omega_0 \beta_i \right)}{2 \left[ (\sigma+ \mu \beta_r - \omega_0 \beta_i)^2 + (\mu \beta_i + \omega_0 \beta_r)^2 \right]},
\end{equation}
see the black curves in Figs. \ref{fig:CGLprofiles}(b)--(d). 
There is good agreement between the onset and buffer curves when \eqref{eq:bufferman} is valid, namely for $\operatorname{Re} \left( \frac{\sigma + (\mu+i \omega_0)(\beta_r+i \beta_i)}{(\mu+i \omega_0)(\beta_r+i \beta_i)} \right) \geq 0$.

The asymptotics are justified by the energy spectral density.
At $\mu = \omega_0$, the energy is concentrated in the low modes, and the $\mathcal{O}(\eps^2 D^2 k^4)$ terms are negligible. 
As $\mu$ increases, the spectrum spreads to higher $|k|$. 
Inclusion of the $\mathcal{O}\left( \eps^2 D^2 k^4 \right)$ terms improves upon \eqref{eq:bufferman} for the buffer curve, better approximating the onset curve at the domain edges (Fig. \ref{fig:CGLprofiles}(c) and (d)).
Fourier series on $[-L,L]$ also yield accurate results for the buffer curves.

As $D \to 0$, \eqref{eq:bufferman} reduces to the prediction for the family of $x$-dependent ODEs. 
Moreover, the ODE for $\hat{A}(\mu,k)$ decouples for each $k$. 
Hence, there is an ODE buffer point for each $k$ \cite{Hayes2016}, and the earliest such point occurs for $k=0$. 
Thus, the constant mode is the first to cause the solution to diverge from the QSS, consistent with \eqref{eq:bufferman}.

For initial conditions with $-\omega_0 <\mu_0 < 0$, \eqref{eq:hom} diverges exponentially before \eqref{eq:inhom}. 
Hence, solutions with $\mu_0$ sufficiently close to $\mathcal{H}$ will transition to oscillations relatively early (near $\mu = -\mu_0$) and not experience maximal delay. 
This is known in ODEs as a memory effect \cite{Baer1989,Neishtadt1987,Neishtadt1988}.

\paragraph{Large-amplitude sources}
The sources in \eqref{eq:CGL} have $\mathcal{O}{(\sqrt{\eps})}$-amplitude as a convenience for analysis. 
For $\mathcal{O}(1)$ and even larger sources, such as $\mathcal{O}(\eps^{-1/2})$, the QSS is a nonlinear function of $\Iapp(x)$. 
(In Fig. \ref{fig:CGLprofiles}(d), the red pencil at the nose of the buffer curve corresponds to the large-amplitude QSS.)
Linearization about the QSS shows that DHB persists, and \eqref{eq:bufferman} is unchanged. 
Fig. \ref{fig:CGLprofiles}(d) demonstrates the effectiveness of the asymptotics for a Gaussian with $a=100$ (so the inhomogeneity in \eqref{eq:CGL} is $\mathcal{O}(\eps^{-1/2})$). 

\TV{
\paragraph{Large diffusivities}
The restriction to weak diffusion, $\eps D$, is also a convenience for analysis. 
For $\eps D = \mathcal{O}(1)$ diffusion with $\mathcal{O}(\sqrt{\eps})$-amplitude sources, the buffer curve prediction \eqref{eq:bufferman} agrees with the numerical onset curve.
}

\paragraph{Summary}
We reported on the persistence of folded node canards and of DHB in the presence of diffusion, extending for the first time these ubiquitous ODE phenomena to PDEs.
The $1^s$ canards of the pituitary cell model \eqref{eq:model} persist for $D>0$ both with spatially homogeneous and inhomogeneous applied currents. 
There are large regimes in parameter space in which one finds steady state solutions with $1^1$
bursting, $1^1 1^0$ alternation, and $1^0$ spiking in different regions.
There are also large regimes with $1^s$ canards for different $s$. 
Also, in the transition intervals between regions, the steady states of the PDE pass through maximal canards.  
Finally, as $D \to 0$, the canards and maximal canards found in the PDE converge to their ODE counterparts.

The spatially-extended pituitary cell model \eqref{eq:model} exhibits delayed loss of stability due to slow passage through HBs (Fig.~\ref{fig:pitslowpass}).
Localized currents lead to spatially inhomogeneous DHB and onset of oscillations, affecting neuronal diagnostics (frequency and amplitude).

To analyse the spatio-temporal DHB, we studied the Complex Ginzburg-Landau PDE \eqref{eq:CGL} with slowly-varying linear growth rate. 
A buffer curve formula that predicts where solutions transition from QSS to oscillatory states was derived.  
In CGL, DHB persists with (up to) $\mathcal{O}(\eps^{-1/2})$ inhomogeneities, and with $\eps D = \mathcal{O}(1)$. 
We found similar results in the Hodgkin-Huxley, Brusselator, FitzHugh-Nagumo, Hindmarsh-Rose, and Morris-Lecar PDEs with source terms, indicating that the new PDE DHB phenomenon is widespread.

In chemical systems (Brusselator and BZ), spatially-localized light sources are sometimes employed. 
DHB in these systems might be studied using the approach developed here. 
Also, the DHB results suggest a new control mechanism in these PDEs.
Rigorous analyses of spatio-temporal canards and DHB in PDEs are in progress.

\begin{acknowledgments}
This research was partially supported by NSF-DMS 1616064. 
We thank R. Bertram, A. Doelman, R. Goh, M. A. Kramer, B. Sandstede, C. E. Wayne, and the referees.
\end{acknowledgments}

\bibliography{Draft}

\end{document}